\documentclass{amsart}

\title[Proper Holomorphic Embeddings]{Proper Holomorphic Embeddings of Finitely and Some Infinitely Connected subsets of $\CC$ into $\CC^2$}
\author{Erlend Forn\ae ss Wold}
\date{January 20, 2005}

\usepackage{amscd}

\newtheorem{theorem}{Theorem}
\newtheorem{lemma}{Lemma}
\newtheorem{proposition}{Proposition}

\theoremstyle{definition}
\newtheorem{defin}{Definition}

\theoremstyle{remark}
\newtheorem{remark}{Remark}

\newcommand{\NN}{\mathbb{N}}
\newcommand{\RR}{\mathbb{R}}
\newcommand{\CC}{\mathbb{C}}
\newcommand{\PP}{\mathbb{P}}

\def\a{{\alpha}}
\def\e{{\epsilon}}
\def\d{{\delta}}
\def\l{{\lambda}}
\def\cO{{\mathcal O}}
\def\r{{\rho}}
\def\s{{\sigma}}

\begin{document}

\begin{abstract}
We show that any finitely connected domain $U\subset\CC$  can be
properly embedded into $\CC^2$.  For some sequences $\{p_j\}\subset
U$, $U\setminus\{p_j\}$   can also be properly embedded into
$\CC^2$.
\end{abstract}

\maketitle

\section{Introduction and Main Result}

\

Let $\mathcal R$  be a noncompact Riemann surface and let
$\phi\colon \mathcal R\rightarrow\CC^2$ be a proper holomorphic
immersion that is 1-1. In that case we say that $\phi$ embeds
$\mathcal R$ \emph{properly} into $\CC^2$.  It is known that any
k-dimensional Stein manifold embeds properly into $\CC^{2k+1}$
[Re, Na, Bi], so in particular any noncompact Riemann surface
embeds properly into $\CC^3$.  It is however an open question
whether any noncompact Riemann surface embeds properly into
$\CC^2$ (it is known that not all compact Riemann surfaces embeds
properly into $\CC\PP^2$, although they do in $\CC\PP^3$ [GH]).
Not much is known even for planar domains. Known results are: The
unit disc by Kasahara and Nishino [KN, St], the annulus by Laufer
[La], the punctured disc by Alexander[Al], and the most general
result so far, due to Globevnik and Stens\o nes: Any finitely
connected bounded domain without isolated points in the boundary.
\u{C}erne and Forstneri\u{c} have some results regarding bordered
Riemann surfaces [\u{C}F]. We prove the following theorem:

\label{main}\begin{theorem} Any finitely connected domain
$U\subset\CC$ can be properly embedded into $\CC^2$.  Moreover, let
$\{p_j\}\subset U$ be a sequence converging to a point $p$ in the
boundary (we allow $p=\infty$), and assume that $\{p_j\}$  is
regular for $U$.  Then $U\setminus\{p_j\}$ can be embedded properly
into $\CC^2$.
\end{theorem}

\

The author would like to thank the referee for many useful comments
and suggestions.

\

\section{Notation And Preliminaries}

\

Throughout this paper we will use the following notation:  For a
real number $R>0$, $B_R$ will denote the open $R$-ball centered at
the origin in $\CC^2$.  We let $\triangle_R$  denote the open
$R$-disc centered at the origin in $\CC$. If there is no subscript
$R$, then we are referring to the unit ball or the unit disc
respectively. We will let $\pi_i$  denote the projection on the i-th
coordinate axis. \

\begin{defin}
Let $U\subset\CC$  be a connected open set such that
$U=\CC\setminus\cup_{i=1}^n K_i$  where the $K_i$'s  are pairwise
disjoint closed connected subset of $\CC$.  We then say that $U$ is
n-connected.
\end{defin}

\begin{defin}
Let $U\subset\CC$  be a domain, and let $\{p_j\}$  be a sequence in
$U$  converging to a boundary point $p$.  We say that $\{p_j\}$  is
\emph{regular}  for $U$  if there exists a continuous curve
$\gamma:[0,1]\rightarrow\overline U$  such that
$\gamma([0,1))\subset U$, with $\gamma(1)=p$, and
$\{p_j\}\subset\gamma$.
\end{defin}

\label{stan}\begin{defin} Let $L=\{\zeta=x+iy\in\CC; x\leq -1,
y=0\}$, let $\Gamma=\{l_j:[0,1]\rightarrow\CC;j=1,...m\}$ be a
collection of smooth disjoint curves without self intersections, and
let $\{p_j\}$ be a discrete set in $\CC$. Assume that these sets are
pairwise disjoint.  Then we will call the domain
$S=\CC\setminus(L\cup\Gamma\cup\{p_i\})$  a \emph{standard} domain.
We will also allow a standard domain to lack some of the above
components in its complement.
\end{defin}

\begin{defin}
Let $Aut_p(\CC^k)$  denote the group of holomorphic automorphisms of
$\CC^k$  fixing the point $p\in\CC^k$.  If all the eigenvalues
$\l_i$ of $\mathrm{d}F(p)$ satisfy $|\l_i|<1$  we say that $F$  is
attracting at $p$.
\end{defin}
\begin{defin}
Let $\{F_j\}\subset Aut_p(\CC^k)$.  We let $F(j)$  denote the
composition map $F_j\circ\cdot\cdot\cdot\circ F_1$, and we define
the basin of attraction of $p$ by
$$
\Omega^p_{\{F_j\}}=\{z\in\CC^k;lim_{j\rightarrow\infty} F(j)(z)=p\}.
$$
\end{defin}

The construction of certain Fatou-Bieberbach domains will be an
integral part of our proof of Theorem \ref{main}, and we will use
the following theorem and proposition:

\label{square}\begin{theorem}\rm{[Wo]} \ Let $0<s<r<1$  such that
$r^2<s$, let $\d>0$, and let $\{F_j\}\subset Aut_0(\CC^k)$  with
$s\|z\|\leq\|F_j(z)\|\leq r\|z\|$ for all $z\in B_\d$, and for all
$j\in\NN$. Then there exists a biholomorphic map
$$
\Phi\colon\Omega^0_{\{F_j\}}\rightarrow\Phi(\Omega^0_{\{F_j\}})=\CC^k.
$$
\end{theorem}
\label{convex}\begin{proposition}\rm[Wo] \ Let $K\subset\CC^k$  be
polynomially convex, let $V\subset\CC^k$  be a closed subvariety,
and let $K'\subset V$ be compact set such that $K\cap V\subset K'$.
Then $\widehat{K\cup K'}_{\cO(\CC^k)}=K\cup\widehat
K'_{\cO(V)}=K\cup\widehat K'_{\mathcal O(\CC^k)}$.
\end{proposition}

\

\section{A Classification of Some Infinitely Connected Domains in
$\CC$}

\

Recall how one can use the existence of a Runge Fatou-Bieberbach
domain $\Omega$  together with The Riemann Mapping Theorem to
embed the unit disc in $\CC$  properly into $\CC^2$.  We may
assume that the intersection between $\Omega$  and the $z$-plane
is not the whole plane, and it follows from the Runge property
that all the connected components of this intersection have to be
simply connected.  Let $U$  be such a component.  Now the Riemann
Mapping Theorem tells us that there is a biholomorphic map $\phi$
mapping $\triangle$  onto $U$.  So if
$\psi\colon\Omega\rightarrow\CC^2$  is the associated
Fatou-Bieberbach map, then $\psi\circ\phi$  will map $\triangle$
properly into $\CC^2$.
\

Now, this method fails if one wants to embed something that is not
simply connected.  This is because intersections between Runge
Fatou-Bieberbach domains and embedded complex curves are simply
connected (it is an open question wether all Fatou-Bieberbach
domains are Runge). Moreover, two multiply connected domains are not
automatically biholomorphically equivalent.  Our approach will be
the following: First we map the domain into $\CC^2$ such that the
image $V$  is Runge, then we construct a Fatou-Bieberbach domain
$\Omega$  such that $V\subset\Omega$  and $\partial
V\subset\partial\Omega$.\

We will prove that any standard domain can be embedded properly
into $\CC^2$, and we begin by investigating which subsets of $\CC$
are biholomorphically equivalent to a standard domain. \

Define the following map $\mu\colon\triangle\rightarrow\CC\setminus
L$:
$$
\mu(\zeta)=(\frac{\zeta+1}{\zeta-1})^2 - 1
$$
For a sequence $p_j$  in $\triangle$  such that $\lim p_j=1$  we
have that $\lim |\mu(p_j)|=\infty$.  In particular, both
$\mu(\triangle)$ and $\mu(\triangle\setminus\{p_j\})$  are standard
domains. \

\label{remark}\begin{remark} The Riemann Mapping Theorem states that
for any simply connected domain $U\subset\CC$ which is not the whole
of $\CC$, there exists a biholomorphism $\phi\colon
U\rightarrow\triangle$ that is onto. Let $p\in\partial U$.  If
$\{p_j\}\subset U$  is a sequence converging to $p$ and if $\{p_j\}$
is regular for $U$, then we can assume that
$lim_{j\rightarrow\infty}\phi(p_j)=1$.
\end{remark}
\

We will use the following theorem from [Go] to show that all the
domains in Theorem \ref{main}  are in fact standard domains. \

\label{hilbert}\begin{theorem}(Hilbert) Every n-connected domain
in the $z$-plane can be mapped univalently onto the $\zeta$-plane
with $n$  parallel finite cuts of inclination $\Theta$  with the
real axis in such a way that a given point $z=a$ is mapped into
$\zeta=\infty$, and the expansion of the mapping function about
$z=a$  has the form
$$
\frac{1}{z-a}+\a_1(z-a)+\cdot\cdot\cdot \  \ \mathrm{or} \ \
z+\frac{\a_1}{z}+\cdot\cdot\cdot
$$
according as a is finite or not.  Some of the cuts referred to may
consist of single points.
\end{theorem}

\

\label{class}\begin{proposition} Let $U\subset\CC$  be n-connected,
and let $\{p_j\}\subset U$  be a sequence converging to a point $p$
in the boundary.  Assume that $\{p_j\}$  is regular for U. Then
$\tilde U=U\setminus\{p_j\}$ is biholomorphically equivalent to a
standard domain $S$.
\end{proposition}
\begin{proof}
Write $U=\CC\setminus\cup_{i=1}^n K_i$, and write $K=\cup_{i=1}^n
K_i$.  In the case of $K$  being unbounded, we will assume that $K$
has only got one unbounded component. It will be clear that the
proof will work also if we have several unbounded components.  We
will have to look at the different possibilities for the limit point
$p$ of the sequence $\{p_j\}$. \

Case 1 - K is bounded, and $p=\infty$:  By Theorem \ref{hilbert}
there is a biholomorphism $\phi$  mapping $U$ onto
 $\CC$  minus a finite number of cuts.  The domain
 $\phi(U)\setminus\{\phi(p_i)\}$  is a standard domain.
\

Case 2 - K is bounded, and $p=K_k$:  By Theorem \ref{hilbert} there
is a biholomorphism $\phi$  mapping $U\cup\{p\}$  onto a domain with
$n-1$  cuts.  Define $\varphi(\zeta)=\frac{1}{\zeta-\phi(p)}$, and
$\varphi\circ\phi$ maps $\tilde U$  onto a standard domain \

Case 3 - K is bounded, and $p\in K_k$:  By Theorem \ref{hilbert}
there is a biholomorphism $\phi$  mapping $U\cup K_k$ onto a domain
with $n-1$  cuts. Define $\varphi(\zeta)=\frac{1}{\zeta-\phi(p)}$.
We have that $\varphi$ maps $\CC\setminus\phi(K_k)$  onto a simply
connected domain $W$ take away zero.  There is a biholomorphism
$f\colon W\rightarrow\triangle$, such that
$\lim_{j\rightarrow\infty} f(\varphi(\phi(p_j)))=1$. So the map
$\mu\circ f\circ\varphi\circ\phi$  maps $\tilde U$  onto a standard
domain. \

Case 4 - K is unbounded, and $p=\infty$:  Let $K_k$ be the unbounded
component of $K$.  By Theorem \ref{hilbert} there is a
biholomorphism $\phi$  mapping $U\cup K_k$ onto a domain with $n-1$
cuts.  Let $V=\overline{\phi(U\cup K_k)}\setminus\phi(K_k)$. $V$ is
simply connected, so there is a biholomorphism $f\colon
V\rightarrow\triangle$ such that
$\lim_{j\rightarrow\infty}f(\phi(p_j))=1$, so the map $\mu\circ
f\circ\phi$  maps $\tilde U$ onto a standard domain. \

Case 5 - K is unbounded, and $p=K_k$:  Define
$\varphi(\zeta)=\frac{1}{\zeta-p}$.  The map $\varphi$  maps $\tilde
U$ onto a domain that falls under Case 1. \

Case 6 - K is unbounded, and $p\in K_k$:  The map
$\varphi(\zeta)=\frac{1}{\zeta-p}$ maps $\tilde U$ onto a domain
falling under Case 4.
\end{proof}

\

\section{Proper Holomorphic Embeddings}

\

In the proof of the following lemma, we use an idea from the proof
of Lemma 2.2 in [BF].\

\label{move}\begin{lemma} Let $K\subset\CC^2$  be a polynomially
convex compact set, let $\e>0$, and let
$\Gamma=\{\gamma_j(t);j=1,..m, t\in [0,\infty)\}$ be a collection of
disjoint smooth curves in $\CC^2\setminus K$  without
self-intersection, such that
$lim_{t\rightarrow\infty}|\pi_1(\gamma_j(t))|=\infty$ for all $j$.
Assume that there exists an $M\in\RR$  such that
$\CC\setminus(\overline\triangle_R\cup\pi_1(\Gamma))$ does not
contain any relatively compact components for $R\geq M$.  Let $p\in
K$. Then for any $R\in\RR$  there exists an automorphism $\phi\in
Aut(\CC^2)$ such that the following is satisfied:

\

(i) $\|\phi(x)-x\|<\e$ for all $x\in K$, \

(ii) $\phi(\Gamma)\subset\CC^2\setminus B_R$,\

(iii) $\phi(p)=p$.
\end{lemma}
\begin{proof}
Let us denote the coordinates on $\CC^2$  by $x=(z,w)$, and the
curves by $\gamma_i(t)=(z_i(t),w_i(t))$ (hence
$z_i(t)=\pi_1(\gamma_i(t))$).  Choose a larger polynomially convex
compact set $K'$  containing an $\e$-neighborhood of K.  We may
assume that $R>M$  and that
$K'\subset\overline\triangle_R\times\CC$.  Let
$\Gamma'=\Gamma\cap(\overline\triangle_R\times\CC)$. \

By Theorem 2.1 in [FL] there exists a $\varphi\in Aut(\CC^2)$ such
that \

(a) $\|\varphi(x)-x\|<\frac{\e}{2}$ for all $x\in K'$,\

(b) $\varphi(\Gamma')\subset\CC^2\setminus\overline B_R$,\

(c) $\varphi(p)=p$.\

To see this, define an isotopy of diffeomorphisms removing $\Gamma'$
from $\overline B_R$.  Since the union of a polynomially convex
compact set and finitely many disjoint smooth compact curves in its
complement is polynomially convex [Sb], the theorem applies. Lastly,
a small translation gives us (c).\

Fix such $\varphi$  and set
$$
\Gamma_R=\{x\in\Gamma;\varphi(x)\in\overline
B_R\}=\Gamma\cap\varphi^{-1}(\overline B_R).
$$
By the construction the complement of
$\overline\triangle_R\cup\pi_1(\Gamma)$  does not contain any
bounded components and
$\pi_1(\Gamma_R)\subset\pi_1(\Gamma)\setminus\overline\triangle_R$.
\

For $i=1,...,m$  let $t_0^i\in\RR^+$  be such that
$z_i(t_0^i)\in\partial\triangle_R$  and
$z_i(t)\in\CC\setminus\overline\triangle_R$  for $t>t_0^i$.  Define
$$
L_i=\{z_i(t_0^i)\}\times\CC,
$$
and let $K_i$ be the intersection
$$
K_i=L_i\cap\varphi^{-1}(\overline B_R).
$$
Since $\varphi^{-1}(\overline B_R)$ is polynomially convex, the
complement of these sets in $L_i$  is connected.  So for any
$T\in\RR$, for each $i$ there is a continuous curve $c_i\colon
[0,1]\rightarrow\CC$ such that $c_i(0)=0$, and such that
$$
\tilde c_i(t)=(z_i(t_0^i),w_i(t_0^i)+c_i(t))
$$
is a curve in $L_i\setminus K_i$  with $|w_i(t_0^i)+c_i(1)|>T$.
There is a neighborhood $U_i$ of $\tilde c_i$  in $\CC^2$  such that
$U_i\subset\subset\CC^2\setminus\varphi^{-1}(\overline B_R)$. For
any $\d>0$  we may now define curves
$$
l_i(t)=(z_i(t_0^i+\d t),w_i(t_0^i+\d t)+c_i(t)).
$$
If we let $\d$  be small enough, the entire curve $l_i$  will be
contained in $U_i$,  and $\pi_2(l_i(1))=T_i$  satisfies $|T_i|> T$.
\

Define the following function $f$  on a subset of $\CC$: \

(i) $f\equiv 0$  on $\overline\triangle_R$,\

(ii) $f(z_i(t_0^i+t))=c_i(\frac{t}{\d})$  for $t\in(0,\d)$,\

(iii) $f(z_i(t_0^i+t))=c_i(1)$  for $t\geq\d$.\

Then $f$  is continuous on
$S=\overline\triangle_R\cup\pi_1(\Gamma)$, and holomorphic on
$\triangle_R$.  For any $C,\r>0$, by Mergelyan's Theorem [Ru] there
exists a holomorphic function $g\in\cO(\CC)$ such that
$\|g-f\|_{S\cap\triangle_C}<\r$.  We may also assume that $g(0)=0$.
Define an automorphism
$$
\psi(z,w)=(z,w+g(z)).
$$
If $\r$ is chosen small enough, and if $T$ and $C$ are chosen big
enough, then each $\psi(\gamma_i)$
 is close to $l_i$  over $z_i(t_0^i+t)$  for $0\leq t\leq\d$, and
$\psi\mid_{\overline\triangle_R\times\CC}\approx id$, such that
$\psi(\Gamma)\cap\varphi^{-1}(\overline B_R)=\emptyset$.  If we also
have that $\r<\frac{\e}{2}$  then $\phi=\varphi\circ\psi$ satisfies
the claims of the lemma.
\end{proof}
\begin{remark}
It is clear from the proof that the corresponding formulation of
Lemma 1 for $\CC^k$  ($k\geq 2$) is also true.
\end{remark}
\label{direction}\begin{lemma} Let
$\gamma\colon[0,1]\rightarrow\CC$ be a $\mathcal C^2$ curve, and
let $\e(\zeta)$  be holomorphic on an open set $U$
 containing $\gamma$. Let $a\in\CC$, and define the following function
$\varphi\colon U\rightarrow\CC^2$:
$$
\varphi(\zeta)=(\zeta,\frac{a}{\zeta-\gamma(0)}+\e(\zeta))
$$
Let $\tilde\gamma(t)=\varphi(\gamma(t))$  for $t\in(0,1]$.  Let
$\pi(t)$  be the projection of the curve $\tilde\gamma$  on the
complex line $L=\{z=w\}$.  Let $c$  be the complex number
corresponding to $\gamma'(0)$, and let
$l(t)=\frac{\sqrt2}{2}(\gamma(0)+\e(\gamma(0))+\frac{a}{ct})$. There
is a $\d>0$ and an $R\in\RR$ such that the following hold: (i)
$|\pi(t)-l(t)|<R$ for $t<\d$, and (ii) $|\pi(t)|$ is strictly
decreasing for $t<\d$.
\end{lemma}
\begin{proof}
$\pi(t)$  is given by
$$
\pi(t)=\frac{\sqrt2}{2}(\gamma(t)+\e(\gamma(t))+\frac{a}{\gamma(t)-\gamma(0)}).
$$
We have to show that the last term gets "close" to
$\frac{\sqrt2}{2}\frac{a}{ct}$ as $t$ gets small. Now
$\gamma(t)=\gamma(0)+ct+h(t)$  where $h(t)=O(\|t\|^2)$, so we have
that
$$
\|\frac{a}{\gamma(t)-\gamma(0)}-\frac{a}{ct}\|=\|\frac{a}{ct+h(t)}-\frac{a}{ct}\|=\|\frac{ah(t)}{c^2t^2(1+\frac{h(t)}{ct})}\|.
$$
Since there is a constant $C$  such that $lim_{t\rightarrow
0}\frac{h(t)}{t^2}=C$, (i) follows. \

Write $x(t)=\frac{2(\gamma(t)-\gamma(0))}{\sqrt2a}$, and
$y(t)=\frac{\sqrt2}{2}(\gamma(t)+\e(\gamma(t)))$.  Now,
$\pi(t)=\frac{1}{x(t)} + y(t)$, so $|\pi(t)|^2$  is given by
$$
m(t)=|\pi(t)|^2=\frac{1}{x(t)\cdot\overline{x(t)}}+2Re(\frac{\overline{
y(t)}}{x(t)}) + y(t)\cdot\overline{y(t)}.
$$
This can be rewritten as
$$
m(t)=\frac{1}{t^{2k}\cdot g(t)} +\frac{h(t)}{t^k\cdot g(t)} +
v(t),
$$
where $h,g,v$  are differentiable real valued functions,
$k\in\NN^+$, $g(0)>0$.  Differentiating $m(t)$  gives $(ii)$.
\end{proof}

\label{embed}\begin{theorem}

Let $U$  be a standard domain as defined in Definition 3.
 Then $U$  can be embedded properly into $\CC^2$.
\end{theorem}

\begin{proof}
We will prove the theorem in the case that $U$  contains all the
sets mentioned in Definition 3 in its boundary, but it will be clear
that the proof works for all standard domains. Let $W=U\cup
L\cup\{l_j(t);t\in (0,1]\}$, and write $q_j=l_j(0)$. We then have
that $W=\CC\setminus\{p_j,q_j\}$.  We start by mapping $W$ properly
into $\CC^2$. Define a function $f\in\cO(W)$ as follows:
$$
f(\zeta)=\sum_{j=1}^m\frac{a_j}{\zeta-q_j}+\sum_{j=1}^\infty\frac{b_j}{\zeta-p_j}.
$$
The coefficients $\{b_j\}$  are chosen such that this converges
uniformly on compact sets in $W$, and such that $f$  is bounded over
$L$.  We will specify conditions on $a_1,...,a_m$  later on.  Now,
we define a map $\omega\colon W\rightarrow\CC^2$ by:
$$
\omega(\zeta)=(\zeta,f(\zeta)).
$$
It is clear that $\|\omega(\zeta)\|\rightarrow\infty$ as
$\zeta\rightarrow q_j$  for $j=1,...,m$, or $\zeta\rightarrow p_j$
for any $j\in\NN$.  We are going to complete the proof by
constructing a Fatou-Bieberbach domain $\Omega$  that intersects the
manifold $\omega(W)$  exactly in $\omega(U)$.  If
$\psi\colon\Omega\rightarrow\CC^2$  is the corresponding
Fatou-Bieberbach map, the composition $\psi\circ\omega$  will map
$U$ properly into $\CC^2$.

Define curves $\gamma_j\colon (0,1]\rightarrow\CC^2$  by
$\gamma_j(t)=\omega(l_j(t))$. Let $L\colon [0,\infty)\rightarrow\CC$
be the curve $L(t)=-1-t$, and let $\gamma_{m+1}(t)=\omega(L(t))$.
Now, since $f$  was chosen to be bounded over $L$, the projection of
$\gamma_{m+1}$  on the plane $\{z=w\}$  is contained in some
$s$-strip around the negative real numbers.  Further, Lemma 2 tells
us that we can choose $a_1,...,a_m$  in the definition of $f$  such
that the projections of the $\gamma_i$'s all point in different
directions.  By the same lemma, none of the projections of the
$\gamma_j$'s are self-intersecting when $t$  is outside of a compact
subset of the open unit interval, so by a change of coordinates, the
conditions on the curves in Lemma 1 are satisfied. \

We will need a polynomially convex compact exhaustion of
$\omega(U)=\omega(W)\setminus\{\gamma_i\}$.  For an $\e>0$, let
$S_\e$ denote
$\{\zeta\in\CC;\mathrm{dist}(\zeta,L\cup\{p_j\}\cup\{l_j\})<\e\}$.
It is clear that if we for an appropriate sequence $\{\e_j\}$
converging to zero, define $C_j=\overline\triangle_j\setminus
S_{\e_j}$,  the sequence $C_1\subset C_2\subset ....$  is a compact
exhaustion of $U$.  The sets $C_i$  are not polynomially convex, but
the sets $K_i=\omega(C_i)$  are. \

We will make one last observation before we start constructing the
Fatou-Bieber\-bach domain:  $\omega(W)$  is a closed submanifold of
$\CC^2$, so by Proposition 1, if $K\subset\CC^2$ is a polynomially
convex compact set such that $K\cap\omega(W)\subset\omega(U)$, then
there exists an $N\in\NN$ such that $K\cup K_i$ is polynomially
convex for all $i\geq N$. \

We may assume that the origin is not contained in any of the
$\gamma_i$'s or in the set $\{p_j\}$, so we can choose a $\r>0$ such
that $\overline B_\r$  does not intersect these sets either. Define
a linear map
$$
A\colon(z,w)\rightarrow(\frac{z}{2},\frac{w}{2}).
$$
Theorem \ref{square} tells us that we can choose a $\d>0$  such that
if we have a sequence of automorphisms $\{\s_k\}\subset
Aut_0(\CC^2)$ such that
\begin{align*}
(*) \ \|\s_k-A\|_{\overline B_\r}<\d
\end{align*}
for all $k\in\NN$, then the basin of attraction to zero of the
sequence $\sigma(k)$  is biholomorphic to $\CC^2$.  What we will do
is to construct a sequence of automorphisms that attracts
$\omega(U)$, but \emph{not} any of the $\gamma_i$'s. \

The sequence of automorphisms will be constructed inductively, and
we make the following induction hypothesis $I_j$:  We have
constructed automorphisms $\{F_1,...,F_j\}\subset Aut_0(\CC^2)$
such that the following is satisfied:

\begin{align}
&\mathrm{Each} \ F_i \ \mathrm{is \ a \ finite \ composition \ of
\ maps} \ \s_k \ \mathrm{satisfying \ (*)},\\
&F(j)(K_j)\subset B_\r,\\
&F(j)(\gamma_i)\subset\CC^2\setminus\overline B_\r \ \mathrm{for} \
i=1,..,m+1.
\end{align}
\

If we define $K_1=\{0\}$  and choose $\r$  small enough, then
$I_1$ is satisfied by letting $F_1=A$.  We will now show how to
construct $F_{j+1}$  so as to ensure $I_{j+1}$. \

$\tilde K=F(j)^{-1}(\overline B_\r)$  is a polynomially convex
compact set satisfying $\tilde K\cap\omega(W)\subset\omega(U)$, so
there exists an $r\geq j+1$  such that $K=\tilde K\cup K_r$ is
polynomially convex. Choose a polynomially convex compact set $K'$
containing a neighborhood of $K$  such that  $K'$ does not intersect
any of the $\gamma_i$'s.\

Choose an $s\in\NN$  such that
$$
(i) \ A^s(F(j)(K'))\subset B_\r,
$$
and then choose an $R\in\RR^+$  such that
$$
(ii) \ A^s(F(j)(x))\notin\overline B_\r \ \mathrm{for \ all} \
\|x\|\geq R.
$$
For any $\e>0$  by Lemma 1 there is a $\phi\in Aut_0(\CC^2)$  such
that $\|\phi(x)-x\|<\e$  for all $x\in K'$, and such that
$$
(iii) \ \phi(\gamma_i)\subset\CC^2\setminus B_R \ \mathrm{for} \
 i=1,...,m.
$$
Define
$$
F_{j+1}=A^s\circ F(j)\circ\phi\circ F(j)^{-1}.
$$
If $\e$  is chosen small enough to ensure that $\phi(K)\subset K'$
then (i) gives us that
$$
F(j+1)(K)=A^s(F(j)(\phi(K)))\subset A^s(F(j)(K'))\subset B_\r,
$$
which ensures (2).  For each $\gamma_i$  by (ii) and (iii) we have
that
$$
F(j+1)(\gamma_i)=A^s(F(j)(\phi(\gamma_i)))\subset
A^s(F(j)(\CC^2\setminus B_R))\subset\CC^2\setminus\overline B_\r,
$$
which ensures (3).  Lastly, the map $A\circ F(j)\circ\phi\circ
F(j)^{-1}$  can be made arbitrarily close to $A$  on $\overline
B_\r$ by choosing $\e$  small enough, which means that $F_{j+1}$
can also be assumed to satisfy (1). \

We have now constructed a sequence of automorphisms $\{F_j\}\subset
Aut_0(\CC^2)$,  and by (1) and Theorem \ref{square}, we have that
$\Omega=\Omega_{\{F_j\}}^0$  is biholomorphic to $\CC^2$.  By (2)
and (3), we have that $\Omega\cap\omega(W)=\omega(U)$. \

Let $\psi$  be the Fatou-Bieberbach map
$\psi:\Omega\rightarrow\CC^2$, and define $\Psi=\psi\circ\omega$.
This map is a proper holomorphic embedding of $U$  into $\CC^2$, and
the proof is finished. \
\end{proof}
Now our main theorem follows easily:

\

\emph{Proof of Theorem 1:} \ By Proposition 2 we have that
$U\setminus\{p_j\}$  is biholomorphically equivalent to a standard
domain $S$.  By Theorem \ref{embed}, S can be embedded properly into
$\CC^2$.

\

Our method has also given us a rather interesting result regarding
Fatou-Bieberbach domains, namely that there exists a
Fatou-Bieberbach domain whose intersection with a complex line $\CC$
is exactly $\CC\setminus L$.

\

\

\

\centerline{Bibliography:}

\

[Al] H.Alexander:{\emph{Explicit imbedding of the (punctured) disc
into $\CC^2$}}, Comment.\ Math.Helv. 52(1977), 539-544.\

[AW] H.Alexander and J. Wermer:{\emph{Several Complex Variables and
Banach Algebras}}, Third edition, Springer-Verlag, New York, 1998.\

[BF] G.Buzzard, F.Forstneri\u{c}:{\emph{A Carleman type theorem for
proper holomorphic embeddings}}, Ark. Mat., 35, 157-169, 1997.\

[Bi] E.Bishop:{\emph{Mappings of partially analytic spaces}}, Amer.
J. Math., 83, 209-242, 1961.\

[CF] M.\u{C}erne, F.Forstneri\u{c}:{\emph{Embedding some bordered
Riemann surfaces in the affine plane}}, Math. Res. Lett., 9,
683-696, 2002.\

[FL] F.Forstneri\u{c}, E.L\o w:{\emph{Global holomorphic equivalence
of smooth manifolds in $\CC^k$}}, Indiana Univ.Math.J, 46, 133-153,
1997.\

[GH] P.Griffiths, J.Harris:{\emph{Principles of Algebraic
Geometry}}, John Wiley and Sons, 1978.\

[Go] D.M.Goluzin: {\emph{Geometric Theory of Functions}},
Translations of Mathematical Monographs, Vol.26, American
Mathematical Society, 1969.\

[GS] J.Globevnik, B.Stens\o nes:{\emph{Holomorphic embeddings of
planar domains into $\CC^2$}}, Math.Ann.303, 579-597, 1995.\

[KN] K.Kasahara, T.Nishino: As announced in Math Reviews 38 (1969)
4721.\

[La] H.B.Laufer:{\emph{Imbedding annuli in $\CC^2$}}, J.d'Analyse
Math, 26(1973) 187-215.\

[MNTU] S.Morosawa, Y.Nishimura, M.Taniguchi,
T.Ueda:{\emph{Holomorphic Dynamics}}, Cambridge University Press,
2000.\

[Na] R.Narasimhan:{\emph{Imbedding of holomorphically complete
spaces}}, Amer.J. Math., 82, 917-934, 1960.\

[Re] R.Remmert:{\emph{Sur les espaces analytiques holomorphiquement
séparables et holomorphiqement convexes}}, C.R. Acad. Sci. Paris,
243, 118-121, 1956.\

[Ru] W.Rudin:{\emph{Real and Complex Analysis}}, McGraw-Hill, Inc,
1987.\

[Sb] G.Stolzenberg:{\emph{Uniform approximation on smooth curves}},
Acta Math. 115, 185-198, 1966.\

[St] J.L.Stehlé:{\emph{Plongements du disques dans $\CC^2$}},
Seminaire P.Lelong(Analyse).Lect.Notes in Math. Vol 275,
Springer-Verlag, Berlin and New York, 1970, 119-30.\

[Wo] E.F.Wold:{\emph{Fatou-Bieberbach domains}}, Preprint, to appear
in International Journal of Mathematics.

\

\centerline{-----------------------------------}

\

Erlend Forn\ae ss Wold\

Ph.d. student, University of Oslo.\

Department of Mathematics\

University of Michigan\

525 E University\

2074 East Hall\

Ann Arbor, MI 48109-1109\

USA\

e-mail: erlendfw@math.uio.no

\end{document}